\documentclass[12pt,reqno]{amsart}

\usepackage{amssymb,latexsym}
\usepackage[matrix,arrow]{xy}

\theoremstyle{plain} 
\newtheorem{thm}{Theorem}[section]
\newtheorem{cor}[thm]{Corollary}

\newtheorem{prop}[thm]{Proposition}
\theoremstyle{remark} 
\newtheorem{rem}[thm]{Remark}
\theoremstyle{definition}

\newcommand{\rl}{\mathbb{R}} \newcommand{\cx}{\mathbb{C}} 
\newcommand{\cl}{\mathbb{C} l}  
\newcommand{\so}{\mathrm{SO}}\newcommand{\Id}{\mathrm{Id}} 
\newcommand{\spin}{\mathrm{Spin}} \newcommand{\wih}[1]{\widehat{#1}} 
\newcommand{\wit}[1]{\widetilde{#1}}  
 \newcommand{\gs}{\sigma} \newcommand{\go}{\omega} 
\newcommand{\gl}{\lambda} \newcommand{\gS}{\Sigma} \newcommand{\gp}{\varphi} 
\newcommand{\gn}{\nabla}   
\newcommand{\ovl}[1]{\overline{#1}}  
\newcommand{\demi}{\frac{1}{2}} \newcommand{\gG}{\Gamma} 
\newcommand{\bS}{\mathbb{S}} \newcommand{\gga}{\gamma} 
\newcommand{\quart}{\frac{1}{4}} 
\newcommand{\ls}{\setlength{\baselineskip}{15pt} 
                 \setlength{\parskip}{3mm}}

\title[]{\bf On Eigenvalue Estimates for the Submanifold Dirac Operator}

\author[]{Nicolas Ginoux and Bertrand Morel} 
\address{Institut {\'E}lie Cartan,
Universit{\'e} Henri Poincar{\'e}, Nancy I,
B.P. 239,
 54506 Vand\oe uvre-L{\`e}s-Nancy Cedex, France}
\email{ginoux@iecn.u-nancy.fr,\;morel@iecn.u-nancy.fr}



\begin{document}
\ls 
\begin{abstract}
We give lower bounds for the eigenvalues of the submanifold Dirac
operator in terms of intrinsic and extrinsic curvature expressions. We also show
that the limiting cases give rise to a class of spinor fields generalizing that of Killing spinors. We conclude by translating
these results in terms of intrinsic twisted Dirac operators.
\end{abstract}
\maketitle
\section{Introduction}

The hypersurface Dirac operator has been used in E.~Witten's proof of the positive energy theorem \cite{Wi}. Lower bounds for the eigenvalues of this operator are given in \cite{Zhang1}, \cite{Zhang2}, \cite{HZ1}. More recently (see \cite{HZM1}), under appropriate boundary conditions, the spectrum of the Dirac operator on compact manifolds with boundary has been examined. In this approach, the hypersurface Dirac operator appears naturally as a boundary term. In \cite{Mor1}, the hypersurface Dirac operator is seen as an intrinsic Dirac-Schr\"odinger operator, with potential given by the mean curvature. 

In this paper, we generalize the results for hypersurfaces obtained in \cite{Mor1} to the case of submanifolds. We start by restricting the spinor bundle of a Riemannian spin manifold to a spin submanifold endowed with the induced metric. We then relate this bundle to the twisted spinor bundle on the submanifold. For convenience, we adapt the algebraic identifications of the spinor spaces and Clifford multiplications given in \cite{Bar1}.

Defining appropriate Dirac operators and relating them with the help of the spinorial Gauss formula, the 
submanifold Dirac operator $D_H$ can be seen as a generalization of the hypersurface Dirac 
operator. As in \cite{HZ2}, we get lower bounds for the eigenvalues of 
$D_H$ in terms of the norm of the mean curvature vector, the Energy-Momentum tensor associated with an eigenspinor, and an adapted conformal change of the metric. However, as a consequence of our definitions, the established estimates hold for all codimensions.

Lower bounds for the eigenvalues of $D_H$ involve the scalar curvature of the submanifold as well 
as a normal curvature term. Assuming that the normal curvature vanishes, the Dirac 
operator on the submanifold twisted with the normal spinor bundle plays the same role as the intrinsic Dirac operator of the hypersurface. 

As in \cite{Mor1}, our identifications allow to discuss the limiting cases in terms of special sections of the spinor bundle. It turns out that these particular spinor fields generalize the notion of Killing spinors to the spinor bundle of the submanifold twisted with the normal spinor bundle.

We finally establish intrinsic analogous results for a modified 
twisted Dirac operator, by considering any auxiliary vector bundle attached to a manifold 
instead of the normal bundle of a submanifold.

We would like to thank Oussama Hijazi for his support during the preparation of this paper.

\section{Dirac operators on submanifolds} 

\subsection{Algebraic Preliminaries}

In this section, we adapt algebraic material developped by C. B\"ar in \cite{Bar1}. Basic facts concerning spinor representations can be found in classical books (see \cite{Friedlivre},\cite{LM},\cite{BFGK} or \cite{BHMM}).

Let $m$ and $n$ be two integers, we start by constructing an irreducible 
representation of the complex Clifford algebra $\cl_{m+n}$ from irreducible 
representations $\rho_n$ and $\rho_m$ of $\cl_n$ and $\cl_m$ respectively. Let 
$\gS_p$ be the space of complex spinors for the representation $\rho_p\,$. Recall 
that if $p$ is even, $\rho_p$ is unique up to an isomorphism, and if $p$ is odd, 
there are two inequivalent irreducible representations of $\cl_p\,$; in this 
case, $(\rho_p^j\,,\gS_p^j)\,$, $j=0,1$, denotes the representation which sends the complex 
volume form to $(-1)^j\,\mathrm{Id}_{\gS_p^j}$. So we have to consider four cases 
according to the parity of $m$ and $n$. 

{\bf First case:} Assume that $n$ and $m$ are even. Define \begin{align*}
 \gga\; 
:\rl^m\oplus\rl^n&\longrightarrow \mathrm{End}_\cx(\gS_m\otimes\gS_n)\\
 (v,w)&\longmapsto \rho_m(v)\otimes(\mathrm{Id}_{\gS_n^+}-\mathrm{Id}_{\gS_n^-})+\mathrm{Id}_{\gS_m}\otimes\rho_n(w), 
\nonumber\end{align*} where $\gS_n^\pm$ is the $\pm 1$-eigenspace for the action of the 
complex volume form $\go_n$ of $\cl_n$. Recall that
$\go_n=i^{[\frac{n+1}{2}]}e_1\cdot\ldots\cdot e_n$, where
$(e_1,\ldots,e_n)$ stands for any positively oriented orthonormal
basis of $\rl^n$ and `$\cdot$' denotes the Clifford multiplication in $\cl_n$.
Then, for $\gs\in\gS_m$, $\theta\in\gS_n$, for any vectors $v\in\rl^m$ and $w\in\rl^n$, we have:
 \begin{eqnarray*}\gga(v+w)^2(\gs\otimes\theta)&=&\rho_m(v)^2\gs\otimes(\theta^++\theta^-)+\rho_m(v)\gs\otimes (\rho_n(w)\theta^--\rho_n(w)\theta^+)\\ 
 && + \rho_m(v)\gs\otimes (\rho_n(w)\theta^+-\rho_n(w)\theta^-)+\gs\otimes \rho_n(w)^2\theta\\
&=& -(|v|^2+|w|^2)\gs\otimes\theta.
\end{eqnarray*}  
Therefore, since 
$\gga(v+w)^2=-|v+w|^2\,\mathrm{Id}$, the map $\gga$ induces a non trivial complex 
representation of $\cl_{m+n}$ of dimension $2^{\frac{m+n}{2}}$ and so $\gga$ is 
equivalent to $\rho_{m+n}$. 

With respect to the inclusions of $\cl_m$ and $\cl_n$ in $\cl_{m+n}$ corresponding to
\begin{align*}
\rl^m\;&\longrightarrow\;\rl^{m+n}=\rl^m\oplus\rl^n\qquad\mathrm{and}\qquad\rl^n\;\longrightarrow\;\rl^{m+n}=\rl^m\oplus\rl^n\\
v\;\;&\longmapsto\;(v,0)\quad\qquad\qquad\qquad\qquad\qquad\;\; w\;\longmapsto\;(0,w),\end{align*}
we can write
\begin{eqnarray}\label{omega}
\go_{m+n}&=& i^{[\frac{m+n+1}{2}]}e_1\cdot\ldots\cdot e_{m+n}\\
&=&i^{[\frac{m+1}{2}]}i^{[\frac{n+1}{2}]}e_1\cdot\ldots\cdot e_m\cdot e_{m+1}\cdot\ldots\cdot e_{m+n}\nonumber\\&=&\go_m\cdot\go_n\, .\nonumber\end{eqnarray}
On the other hand, if $\gs\in\gS_m$ and $\theta\in\gS_n$, then for all $v\in \rl^m$,  \begin{equation}\label{vgo1}
\gga(v\cdot\go_n)(\gs\otimes\theta)=\rho_m(v)\gs\otimes\theta \,.\end{equation}
Therefore, since $m$ is even, we have \begin{equation*}\gga(\go_{m+n})(\gs\otimes\theta)=\rho_m(\go_m)\gs\otimes\rho_n(\go_n)\theta \end{equation*} so that
\begin{align*} \gS_{m+n}^+&=\gS_m^+\otimes\gS_n^+\oplus\gS_m^-\otimes\gS_n^- \\
\gS_{m+n}^-&=\gS_m^+\otimes\gS_n^-\oplus\gS_m^-\otimes\gS_n^+ \;. \end{align*}
We can then define $$\gS:=\gS_m\otimes\gS_n=\gS_{m+n}^+\oplus\gS_{m+n}^-.$$
{\bf Second case :} Assume that $m$ is odd and $n$ is even. For $j=0,1$, set
 \begin{align*}
 \gga^j\; 
:\rl^m\oplus\rl^n&\longrightarrow \mathrm{End}_\cx(\gS_m^j\otimes\gS_n)\\ 
 (v,w)&\longmapsto \rho_m^j(v)\otimes(\mathrm{Id}_{\gS_n^+}-\mathrm{Id}_{\gS_n^-})+\mathrm{Id}_{\gS_m^j}\otimes\rho_n(w), 
\nonumber\end{align*} As before,
 the map $\gga^j$ induces a non trivial complex representation of 
$\cl_{m+n}$ of dimension $2^{[\frac{m+n}{2}]}$. Since $\go_{m+n}=\go_m\cdot\go_n$ as in (\ref{omega}), we have
$\gga^j(\go_{m+n})=(-1)^j\,\mathrm{Id},$ and therefore the representations $\gga^j$ and  $\rho_{m+n}^j$ are equivalent.  
Note that \begin{equation}\label{vgo2}
\gga^j(v\cdot\go_n)=\rho_m^j(v)\otimes\Id_{\gS_n}\; , \quad \forall v\in \rl^m \,.\end{equation}
{\bf Third case :} Assume that $m$ is even and $n$ is odd. For $j=0,1$, set 
 \begin{align*}
 \gga^j\; 
:\rl^m\oplus\rl^n&\longrightarrow \mathrm{End}_\cx(\gS_m\otimes\gS_n^j)\\[0.3cm] 
 (v,0)&\longmapsto i \begin{pmatrix}
 0&-\rho_m(v)\\\rho_m(v)&0 \end{pmatrix}\otimes\Id_{\gS_n^j}\nonumber\\[0.3cm](0,w)&\longmapsto \begin{pmatrix}
 \Id_{\gS_m^+}&0\\0&-\Id_{\gS_m^-}
 \end{pmatrix}\otimes\rho_n^j(w), \nonumber
\end{align*}
where the matrices are given with respect to the decomposition 
$\gS_m=\gS_m^+\oplus\gS_m^-$. Once again, $\gga^j$ is an 
irreducible complex representation of $\cl_{m+n}$. As in the previous case, $\go_{m+n}=\go_m\cdot\go_n$ and we see that $\gga^j(\go_{m+n})=(-1)^j\mathrm{Id}\,.$

 So we proved that $\gga^j$ is equivalent to $\rho_{m+n}^j$ and \begin{equation}\label{vgo3}
\gga^j(v\cdot\go_n)=(-1)^j\,
i\,\rho_m^j(v)\otimes\Id_{\gS_n}\; , \quad \forall v\in \rl^m\, .\end{equation}
{\bf Fourth case :} Assume that $m$ and $n$ are odd. Define \begin{align*}
\gS^+&:=\gS_m^0\otimes\gS_n^0\, ,\\\gS^-&:=\gS_m^0\otimes\gS_n^1\, 
,\\\gS&:=\gS^+\oplus\gS^-\, ,
\end{align*}
and
 \begin{align*}
 \gga\; 
:\rl^m\oplus\rl^n&\longrightarrow \mathrm{End}_\cx(\gS)\\[0.3cm] 
 (v,0)&\longmapsto i \begin{pmatrix}
 0&\rho_m^0(v)\otimes\tau^{-1}\\-\rho_m^0(v)\otimes\tau&0 \end{pmatrix}\nonumber\\[0.3cm](0,w)&\longmapsto \begin{pmatrix}
0&-\Id_{\gS_m^0}\otimes \tau^{-1}\circ\rho_n^1(w)\\\Id_{\gS_m^0}\otimes \tau\circ\rho_n^0(w)&0
 \end{pmatrix},\nonumber 
\end{align*}
where $\tau$ is an isomorphism from $\gS_n^0$ to $\gS_n^1$ satisfying 
\begin{equation*}
 \tau\circ\rho_n^0(w)\circ\tau^{-1}=-\rho_n^1(w)\, , \quad \forall w\in\rl^n.
\end{equation*} 
Now, as in previous cases, we have $\gga(v+w)^2=-(|v|^2+|w|^2)\mathrm{Id}_\gS$ for all $v\in\rl^m$ and $w\in\rl^n$. Moreover, since in the case where $m$ and $n$ are odd, $\go_{m+n}=-i\,\go_m\cdot\go_n\,$, we can show that $$\gga(\go_{m+n})=\begin{pmatrix}
\Id_{\gS^+}&0 \\ 0&-\Id_{\gS^-}\end{pmatrix}\;.$$ Therefore we conclude that $\gga$ is equivalent to $\rho_{m+n}$ 
and $\gS_{m+n}^\pm\cong \gS^\pm\,$.

Besides, we have the relation \begin{equation}\label{vgo4} \gga(v\cdot\go_n)=i\begin{pmatrix}\rho^0_m(v)\otimes\Id_{\gS_n^0}&0\\0&-\rho^0_m(v)\otimes\Id_{\gS_n^1}\end{pmatrix},\quad \forall v\in \rl^m\, .\end{equation} 

\subsection{Restriction of Spinors to a Submanifold} 

Let $(\wit{M}^{m+n},g)$ be a Riemannian spin manifold and let $M^m$ be an immersed 
oriented submanifold in $\wit{M}$ with the induced Riemannian structure. Assume that $(M^m,g_{|M})$ is spin. If $NM$ is the normal vector bundle of $M$ in 
$\wit{M}$, then there exists a spin structure on $NM$, denoted by $\spin N$. Let $\spin M\times_M \spin 
N$ be the pull-back of the product fibre bundle $\spin M\times\spin N$ over $M\times M$ by the diagonal map. There exists 
a principal bundle morphism ${\Phi: \spin M\times_M \spin 
N\rightarrow\spin\wit{M}_{|M}}$,  with \begin{equation}\label{action} 
  \Phi((s_M,s_N)(a,a'))= \Phi((s_M,s_N)) (a\cdot a')
\end{equation} for all $(s_M,s_N)$ in $\spin M\times_M \spin N$ and for all 
$(a,a')$ in $\spin(m)\times\spin(n)$, such that the following diagram commutes : $$\xymatrix{\spin M\times_M \spin 
N\ar[dd]\ar[r]^{\qquad\Phi}&\spin \wit{M}_{|M}\ar[dd]\ar[dr]&\\&&M\\ \so 
M\times_M \so N \ar[r]&\so\wit{M}_{\mid M}\ar[ur]&}$$ where the lower horizontal 
arrow is just given by juxtaposition of bases (see \cite{Miln1}). 

Now, let $\bS:=\gS\wit{M}_{\mid M}$, where $\gS\wit{M}$ is the spinor 
bundle of $\wit{M}$ and 
\begin{equation*}
  \gS:=\begin{cases}
 \gS M\otimes\gS N& \mbox{if } n \mbox{ or } m \mbox{ is even,}\\
 \gS M\otimes\gS N\oplus \gS M\otimes\gS N & \mbox{otherwise}.
  \end{cases}
\end{equation*}
 
Recall that there exists a hermitian inner product on $\bS$, denoted by 
$<.\, ,.>$, such that Clifford 
multiplication by a vector of $T\wit{M}_{\mid M}$ is
skew-symmetric. In the following, we write  $(.\, ,.)=\Re e(<.\, , .>)$. 

\subsection{Identification of the Restricted Spinor Bundle} 

{}From the preceding considerations, it is now possible to identify $\bS$ with $\gS$. 
For example, if $m$ and $n$ are even, we have the following isomorphism:
\begin{eqnarray*}
\gS M \otimes \gS N &\longrightarrow & \bS\nonumber\\
([s_M,\gs],[s_N,\eta]) &\longmapsto & [\Phi(s_M,s_N), \gs\otimes\eta]
\end{eqnarray*}
where the last equivalence class is given, for all $(a,a')\in 
\spin(m)\times\spin(n)$, by
\begin{equation*}
   \Big(\Phi((s_M,s_N)(a,a')),\gs\otimes\eta\Big)\sim\Big(\Phi(s_M,s_N),\gga 
   (a\cdot a')(\gs\otimes\eta)\Big),
\end{equation*} with respect to (\ref{action}). From now on, the inverse of this isomorphism will be denoted by 
\begin{equation}\label{etoile}
\psi\in\gG(\bS)\mapsto\psi^\star\in\gG(\gS).
\end{equation}
 With respect to ${<.\, ,.>}$ and the naturally induced hermitian inner product on 
$\gS$, this isomorphism is unitary. This is why both inner 
products will be denoted by the same symbol when using this identification. 

Let $\go_\perp=\go_n$ if $n$ is even, and 
$\go_\perp=-i\go_n$ if $n$ is odd. Recall that in both cases ${\go_\perp^2=(-1)^n}$ (compare with the definition of $\go_\perp$ in \cite{HZ2} and note that it keeps the same properties). From (\ref{vgo1}), (\ref{vgo2}), (\ref{vgo3}) and (\ref{vgo4}), it is easy to see that, 
with respect to the representation $\gga$ defined in Section 2.1, Clifford multiplication by a vector field $X$ tangent to $M$ satisfies
 \begin{equation}\label{propmult} 
   \forall \psi\in\Gamma(\bS),\quad 
  X\underset{M} \cdot\psi^\star=(X\cdot\go_\perp\cdot\psi)^\star.
\end{equation}
 
\subsection{The Gauss Formula and the Submanifold Dirac Operator} 

Fix $p\in M$ and denote by $(e_1,\ldots,e_m,\nu_1,\ldots ,\nu_n)$ a positively 
oriented local orthonormal basis of $T\wit{M}_{\mid M}$ such that 
$(e_1,\ldots,e_m)$ (resp. $(\nu_1,\ldots,\nu_n)$) is a positively
oriented local orthonormal basis of $TM$ (resp. $NM$). If $\wit{\gn}$ denotes the Levi-Civita 
connection of $(\wit{M},g)$, then for all 
$X\in\gG(TM)$, for all $Y\in \gG(NM)$ and for $i=1,\ldots,m$, the Gauss formula can be written as 
\begin{equation}\label{Gauss} 
\wit{\gn}_i (X+Y)=\gn_i(X+Y)+h(e_i,X)-h^*(e_i,Y),
\end{equation} 
where $\gn_i(X+Y)=\gn^M_i X+\gn^N_iY$, and $h^*(e_i,\,.)$ is the transpose of the second fundamental form $h$ 
 viewed as a linear map from $TM$ to $NM$. Here $\wit{\gn}_i$ stands 
for $\wit{\gn}_{e_i}$.

Denote also by $\wit{\gn}$ and $\gn$ the induced spinorial covariant
derivatives on $\gG(\bS)$. Therefore, on $\gG(\bS)$, $\gn=\gn^{\gS M}\otimes\Id + \Id\otimes\gn^{\gS N}$ except for $n$ and $m$ odd where 
$\gn=(\gn^{\gS M}\otimes\Id + \Id\otimes\gn^{\gS N})\oplus 
(\gn^{\gS M}\otimes\Id + \Id\otimes\gn^{\gS N})$. For
$\psi\in\gG(\bS)$, the covariant derivative 
$\gn\psi$ is understood via the relation $(\gn\psi)^\star=\gn\psi^\star$. 

As in \cite{Bar1}, one can deduce from (\ref{Gauss}) the spinorial Gauss 
formula:
\begin{equation}\label{Sgauss}
  \forall \psi \in \gG(\bS),\quad  \wit{\gn}_i\psi=\gn_i\psi+\frac{1}{2}\sum_{j=1}^m e_j\cdot h_{ij}\cdot\psi.
\end{equation}
Now, define the following Dirac operators
\begin{equation*}
\wit{D}=\sum_{i=1}^{m} e_i\cdot\wit{\gn}_i\; ,\quad D=\sum_{i=1}^m e_i\cdot\gn_i\;,\end{equation*}
and, $H=\sum_{i=1}^m h(e_i,e_i)$ denoting the mean curvature vector field,
\begin{equation}\label{dh}
  D_H:=(-1)^n\go_\perp\cdot\wit{D}=(-1)^n\go_\perp\cdot D + \demi H\cdot\go_\perp\cdot\psi 
\end{equation}
since $H\cdot\go_\perp\cdot=(-1)^{n-1}\go_\perp\cdot H\cdot\;$ and $  \wit{D}=D-\demi H\cdot\;$ by (\ref{Sgauss}).

\begin{rem}
Another Dirac operator can be defined by using intrinsic Clifford multiplication and twisting the Dirac 
operator on the submanifold with the normal spinor bundle. This has 
been done by C. B\"{a}r in \cite{Bar1} by setting 
\begin{eqnarray}
D_M^{\gS N}\,:\,\gG(\gS)&\longrightarrow&\gG(\gS)\nonumber\\ \gp&\longmapsto&\begin{cases} \sum_i 
e_i\underset{M}  \cdot\gn_i \gp &\mbox{ if } m \mbox{ or } n \mbox{ is even,}\\ \sum_i e_i\underset{M}  \cdot 
\gn_i \gp \oplus -\sum_i e_i \underset{M} \cdot \gn_i \gp &\mbox{ if } m \mbox{ and } n \mbox{ are 
odd.}\end{cases}\nonumber\end{eqnarray} 

In fact, with the help of (\ref{propmult}) and (\ref{dh}), we can relate $D_H$ and $D_M^{\gS N}$ by
\begin{eqnarray}\label{CompareD}
  (D_H\psi)^\star&=& \Big((-1)^n \go_\perp\cdot D\psi + \demi H\cdot\go_\perp\cdot\psi\Big)^\star\nonumber\\ 
&=& D_M^{\gS N}\psi^\star+\demi (H\cdot\go_\perp  \cdot \psi)^\star ,\quad\forall \psi \in \gG(\bS).
\end{eqnarray}
\end{rem} 

It is known that $D_H$ is formally self-adjoint and that $D_H^2=\wit{D}^{\,*}\wit{D}$, 
where $\wit{D}^{\,*}$ is the formal adjoint of $\wit{D}$ w.r.t. $\int_M(.\, ,.)v_g$ (see 
\cite{HZ2}). 
 
\section{Estimates for the Eigenvalues of the Submanifold Dirac Operator}

\subsection{Basic Estimates}
First, for any spinor field $\psi \in \gG(\bS)$, define the function  
\begin{equation}\label{defrn}
 R^N_\psi:= 2\sum_{i,j=1}^m (e_i\cdot e_j\cdot\mathrm{Id}\otimes \mathrm{R}^N_{e_i,e_j}\psi,\psi/|\psi|^2) 
\end{equation} on $M_\psi:=\{x\in M\;:\;\psi(x)\neq 0\}$, where $\mathrm{R}^N_{e_i,e_j}$ stands for spinorial normal
curvature tensor. We start by giving a proof of the
following result (see \cite{HZ2}):
\begin{thm}[Hijazi-Zhang]\label{thm1}
Let $M^m \subset \wit{M}^{m+n}$ be a compact spin submanifold of a Riemannian spin manifold 
$(\wit{M},g)$. Consider a non-trivial spinor field $\psi \in \gG(\bS)$ such that $D_H\psi=\gl\psi$. Assume that $m\geq 2$ and $$ m 
(R+R^N_\psi)>(m-1) ||H||^2 > 0$$ on $M_\psi$, where $R$ is the scalar curvature of $(M^m,g_{|M})$ and $R^N_\psi$ is given by (\ref{defrn}).
Then one has 
\begin{equation}\label{ethm1}
\gl^2 \geq \frac{1}{4}\; \inf _{M_{\psi}} \left( \sqrt {\frac{m}{m-1}(R+R^N_\psi)} -||H|| \right) ^2\; .
\end{equation}
\end{thm}

\begin{proof}
For any function $q$, nowhere equal to $\frac{1}{m}$, define the modified connection,  
\begin{equation*}
\gn _i ^\gl=\gn _i +\frac{1-q}{2(1-mq)}\, e_i \cdot H \cdot + q\gl e_i \cdot \go_\perp\cdot\,.
\end{equation*} 
Using the Lichnerowicz-Schr\"{o}dinger formula (see \cite{LM}) we have $$(D^2\psi,\psi)=(\gn^*\gn\psi,\psi) + 
\quart(R+R^N_\psi)|\psi|^2,$$ and one can easily compute
\begin{align}\label{presq}
\int_M& |\gn^\gl\psi|^2 v_g=\nonumber\\
&\int_M (1+mq^2-2q)\Big[\gl^2-\frac{1}{4}\Big(\frac{R+R^N_\psi}{(1+mq^2-2q)}
-\frac{(m-1)||H||^2}{(1-mq)^2}\Big)\Big] |\psi |^2 v_g 
\end{align}

Then, assuming $m(R+R^N_\psi)>(m-1) ||H||^2 > 0$ on $M_{\psi}$, we can choose $q$ so that 
\begin{equation}\label{defq}
  (1-m q)^2= \frac{(m-1)||H||}{\sqrt{\frac{m}{m-1}(R+R^N_\psi)} -||H||}\quad\textrm{ on }M_\psi.
\end{equation}
Inserting equation (\ref{defq}) in (\ref{presq}), and since the complement of $M_\psi$ in $M$ is of zero-measure, we conclude by observing that the left 
member of (\ref{presq}) is nonnegative. 
\end{proof}

Let $\kappa_1$ be the lowest eigenvalue of the self-adjoint operator $\mathcal{R}^N$ defined by
\begin{eqnarray}
  \mathcal{R}^N\;:\;\gG(\bS)&\longrightarrow&\gG(\bS)\\
  \gp&\longmapsto& 2\sum_{i,j=1}^m e_i\cdot e_j\cdot\mathrm{Id}\otimes \mathrm{R}^N_{e_i,e_j}\gp\; .
\end{eqnarray}
The hypothesis $m(R+R^N_\psi)>(m-1) ||H||^2 > 0$ in Theorem \ref{thm1} can be strengthened to give

\begin{cor}\label{cor1}
Under the same hypotheses as in Theorem \ref{thm1}, assume that $m\geq 2$ and $$m(R+\kappa_1)>(m-1) ||H||^2 > 0$$ on $M$,
 then 
\begin{equation}\label{ecor1}
\gl^2 \geq \frac{1}{4}\; \inf _M \left( \sqrt {\frac{m}{m-1}(R+\kappa_1)} -||H|| \right) ^2\; .
\end{equation}
\end{cor}

Recall that in the case of hypersurfaces, limiting cases are characterized by the existence of a real Killing spinor on $M$ and the fact that the mean curvature $H$ is constant (see \cite{Zhang1} and \cite{Mor1}). A non-zero section $\psi$ of $\,\bS$ satisfying \[\forall X\in\Gamma(TM),
\quad\gn_X\psi^\star=-\frac{\mu}{m}X\underset{M}\cdot\psi^\star\] for a given real constant $\mu$ will be called a twisted (real) Killing spinor. 

\begin{prop}
If equality holds in (\ref{ecor1}), then $(M^m,g_{|M})$ admits a
twisted Killing spinor and $||H||$ is constant.\end{prop}  
\begin{proof}
Suppose the limiting case holds in (\ref{ecor1}), then
the right hand side has to be constant on $M$, and 
\begin{equation}
\gl^2=\frac{1}{4}(\sqrt{\frac{m}{m-1}(R+\kappa_1)}-||H||)^2\; ,\qquad
\gn^\gl \psi=0,\quad\textrm{ on }M\label{par}.
\end{equation}
Note that equality holds in (\ref{ethm1}) which yields $R^N_{\psi}=\kappa_1\;$. Hence $\psi$ is an eigenspinor for the operator $\mathcal{R}^N$ with eigenvalue $\kappa_1$.
Using (\ref{par}), we can show that $|\psi|$ must be constant on $M$ (therefore, $M_\psi=M$) and compute 
\begin{eqnarray*}
D\psi&=&-\sum_{i=1}^m e_i\cdot\Big (\frac{1-q}{2(1-mq)}\, e_i \cdot H \cdot\psi + q\gl e_i \cdot \go_\perp\cdot\psi\Big )\\
&=&\frac{m(1-q)}{2(1-mq)}H\cdot\psi + m q\gl \go_\perp\cdot\psi
\, .
\end{eqnarray*} Then, by (\ref{dh}) and the fact that
$H\cdot\go_\perp\cdot=(-1)^{n-1}\go_\perp\cdot H\cdot$, 
\begin{eqnarray*}
0&=&\gl\go_\perp\cdot\psi+\frac{H}{2}\cdot\psi-\frac{m(1-q)}{2(1-mq)}H\cdot\psi-m q\gl \go_\perp\cdot\psi\\
0&=&(1-mq)^2\gl\go_\perp\cdot\psi-(m-1)\frac{H}{2}\cdot\psi
\, .
\end{eqnarray*} Since in the equality case, $\sqrt{\frac{m}{m-1}(R+\kappa_1)}-||H||=2|\gl|$, we can deduce the relation:
\[ \go_\perp\cdot\psi=\mathrm{sgn}(\gl)\frac{H}{||H||}\cdot \psi.\]

With respect to the isomorphism ``$\;{}^\star\;$'', we can rewrite equation (\ref{par}) as an intrinsic equation on $\gG(\gS)$: \begin{eqnarray*}
\forall X \in \gG(TM),\quad \gn_X\psi^\star&=&-\frac{\mu}{m}X\underset{M}\cdot\psi^\star\\
\textrm{ with }\quad\mu&=&\frac{\mathrm{sgn(\gl)}}{2}\sqrt{\frac{m}{m-1}(R+\kappa_1)}.
\end{eqnarray*}
Note that if there exists two smooth real
functions $f$ and $\kappa$ on $M$ and a non-zero
section $\psi$ of $\,\bS$ satisfying for all vector field $X$ on $M$ $$\quad\gn_X\psi^\star=-\frac{f}{m}X\underset{M}\cdot\psi^\star\quad \mathrm{and} \quad\mathcal{R}^N\psi=\kappa\psi\;,$$ then, by
computing the action of the curvature tensor on $\psi^\star$, we see that
necessarily 
\begin{equation*}
 \demi\mathrm{Ric}(X)\underset{M}\cdot\psi^\star-\sum_{i=1}^m (e_i\cdot\mathrm{Id}\otimes\mathrm{R}_{X,e_i}^N)\psi^\star=-\frac{1}{m}\mathrm{d}f\underset{M}\cdot X\underset{M}\cdot\psi^\star- \mathrm{d}f(X)\psi^\star+2\frac{m-1}{m^2}f^2 X\underset{M}\cdot\psi^\star 
\end{equation*}
which implies
\[\qquad f^2=\frac{m}{4(m-1)}(R+\kappa)=\mathrm{constant}\;.\]
Moreover, in the equality case, the fact that $f$ is constant implies that $||H||$ is constant.\end{proof}

\begin{rem}
If the normal curvature tensor is zero, then $\mu$ has to be constant and the manifold $M$ must be Einstein with mean curvature
vector being of constant length. Besides, the equality case
corresponds to that of Friedrich's inequality. Therefore $\mu$ is the first
eigenvalue of the Dirac operator $D_M^{\gS N}$.
\end{rem}

\subsection{Estimate Involving the Energy-Momentum Tensor}

If $\psi \in \gG(\bS)$ is a spinor field, we define the Energy-Momentum tensor $Q^{\,\psi}$ 
associated with $\psi$ on $M_\psi$ by 
$$Q^{\,\psi}_{ij}=\frac{1}{2} (e_i \cdot \go_\perp \cdot\gn_j \psi + e_j \cdot\go_\perp 
\cdot\gn_i \psi, \psi/|\psi |^2).$$ 
Note  that 
\begin{equation*}
 Q^{\,\psi}_{ij}=\frac{1}{2} (e_i \underset{M} \cdot\gn_j \psi^\star + e_j  
\underset{M}\cdot\gn_i \psi^\star, \psi^\star/|\psi^\star |^2).
\end{equation*}
Therefore, $Q^{\psi}$ is the intrinsic Energy-Momentum tensor
associated with $\psi^\star$. Observe that this intrinsic Energy-Momentum tensor is the one that appears in the Einstein-Dirac equation (see \cite{FK1}). We prove the following (compare
with \cite{Mor1})
\begin{thm}\label{thm2}
Let $M^m \subset \wit{M}^{m+n}$ be a compact spin submanifold of a Riemannian spin manifold 
$(\wit{M},g)$. Consider a non-trivial spinor field $\psi \in \gG(\bS)$ such that $D_H\psi=\gl\psi$. Assume that $$R+\kappa_1+4|Q^{\,\psi}|^2> ||H||^2 > 0$$ on $M_\psi$.
Then one has 
\begin{equation}\label{ethm2}
\gl^2 \geq \frac{1}{4}\; \inf _{M_{\psi}} \left( \sqrt {R+\kappa_1+4|Q^{\,\psi}|^2} -||H|| \right) ^2. 
\end{equation} 
\end{thm}
\begin{proof}
For any real function $q$ that never vanishes, consider the modified 
covariant derivative defined on $\gG(\bS)$ by 
\begin{eqnarray*}
\gn _i ^Q=\gn _i - \frac{1}{2mq} e_i\cdot H\cdot +(-1)^{n+1}q\gl e_i \cdot \go_\perp \cdot +\,\sum_j 
Q^{\,\psi}_{ij}\;e_j\cdot\go_\perp \cdot \,. 
\end{eqnarray*}
As in the proof of Theorem \ref{thm1}, we compute
\begin{align}
  \label{Qpres}
 \int_M& |\gn^Q \psi |^2v_g=\nonumber\\ &\int_M
(1+mq^2)\Big[\gl^2-\frac{1}{4}\Big(\frac{R+R^N_\psi+4|Q^{\,\psi}| 
^2}{(1+mq^2)}-\frac{||H||^2}{mq^2}\Big )\Big] |\psi |^2 v_g\nonumber\\-&\frac{1}{4}\int_M(1+mq^2)\Big[ \frac{2}{mq(1+mq^2)}( ||H||^2-\frac{<H\cdot\psi,\go_\perp\cdot\psi>^2}{|\psi|^4})\Big] |\psi |^2 v_g 
\end{align}

To finish the proof of Theorem \ref{thm2}, if
$R+\kappa_1+4|Q^{\,\psi}|^2> ||H||^2 >0$, we take 
\begin{equation*}
  q=\sqrt{\frac{||H||}{m(\sqrt{R+\kappa_1+4|Q^{\,\psi}| ^2}-||H||)}}\;,
\end{equation*}
and then observe that by the Cauchy-Schwarz inequality, we have
\begin{equation}
  \label{CS}
  ||H||^2-\frac{<H\cdot\psi,\go_\perp\cdot \psi>^2}{|\psi|^4}\geq 0.
\end{equation}
\end{proof}

Suppose now that equality holds in (\ref{ethm2}). Then 
\begin{equation*}
\gn^Q\psi=0\quad,\quad|\gl|=\frac{1}{2}\;  \left( \sqrt {R+\kappa_1+4|Q^{\,\psi}|^2} -||H|| \right)\quad\mathrm{and}\quad
\mathcal{R}^N\psi=\kappa_1\psi\; .
\end{equation*}
Moreover,
\[ ||H||^2-\frac{<H\cdot\psi,\go_\perp\cdot \psi>^2}{|\psi|^4}= 0,\]
so that, by the equality case in the Cauchy-Schwarz inequality, 
\[\go_\perp\cdot\psi=f H\cdot\psi,\] for some real function
$f$ on $M$. As in the preceding section, and taking into account the identification (\ref{etoile}), we deduce that
$f=\frac{\mathrm{sgn}(\gl)}{||H||},$ and that the section $\psi$ satisfies
\begin{equation}\label{EMspinor}\gn_i\psi^\star=-\sum_{j}Q^\psi_{ij}\,e_{j}\underset{M}\cdot\psi^\star.\end{equation}
Hence, we can say that $\psi$ is a kind of Energy-Momentum spinor (see 
\cite{Mor1}). We will call such a section a twisted EM-spinor. One can give an integrability condition for the existence of twisted EM-spinors, by computing the action of the curvature tensor on $\gG(\bS)$:
\begin{equation*}
(\mathrm{tr}(Q^\psi))^2=\frac{1}{4}(R+R^N_\psi+4|Q^\psi|^2).
\end{equation*}
This implies, with equation (\ref{EMspinor}), that the section $\psi^\star$ is an ``eigenspinor'' for $D_M^{\gS N}$ associated with the function $\pm \demi\sqrt{R+\kappa_1+4|Q^\psi|^2}$.
Note that this function is constant if and only if $||H||$ is constant.

\subsection{Conformal Lower Bounds}

Consider a conformal change of the metric $\ovl{g}=e^{2u}g$ for a real 
function $u$ on $\wit{M}$. Let  
\begin{eqnarray}
  \label{chang}
  \bS\,&\longrightarrow\ovl{\bS}\\\psi\,&\longmapsto \ovl{\psi}\nonumber\end{eqnarray}
be the induced isometry between the two corresponding spinor bundles. Recall that if
$\gp$, $\psi$ are two sections of $\bS$, and $Z$ any vector field on $\wit{M}$, we have 
$$(\gp,\psi)=(\ovl{\gp},\ovl{\psi})_{\ovl{g}}\quad \mathrm{and}\quad 
\ovl{Z}\;\bar{\cdot}\;\ovl{\psi}=\ovl{Z\cdot\psi}, $$  
where $\ovl{Z}=e^{-u}Z$. We will also denote by $\ovl{g}=(e^{2u}g)_{\mid M}$ the restriction of $\ovl{g}$ 
to $M$. 

Note that this isomorphism commutes with the isomorphism
``$\;{}^\star\;$'' given by (\ref{etoile}). By conformal covariance of the Dirac operator, for $\psi \in 
 \gG(\bS)$, we have,
\begin{equation}\label{covconf}
\overline D \Big(\; e ^{-\frac{(m-1)}{2} u} \;\overline{\psi}\; \Big)=e 
^{-\frac{(m+1)}{2} u} \;\overline {D \psi},  
\end{equation}
where $\overline D$ stands for the Dirac operator w.r.t. to $\ovl{g}$. On the 
other hand, the corresponding mean curvature vector field is given by 
\begin{equation}\label{Hbar}
\wit{H} = e ^{-2u} \,\Big( H- m\; \mathrm{grad}^N\,u \Big). 
\end{equation}
Now, assume that $\mathrm{grad}^N\,u=0$. If $\overline{D}_{H}$ stands for the submanifold Dirac 
operator w.r.t. to $\ovl{g}$, equations (\ref{covconf}) and (\ref{Hbar}) imply 
that $\overline{D}_{H}$ is a conformally covariant operator, i.e. 
\begin{equation}\label{DHbar}
\ovl{D}_{H}\; \Big(\; e ^{-\frac{(m-1)}{2} u}\; \ovl{\psi}\; \Big) = e 
^{-\frac{(m+1)}{2} u }\; \Big(\; {\overline {D _H \psi}} \;\Big). 
\end{equation} for any section $\psi$ of $\bS$. 

>From now on, we will only consider regular conformal changes of the metric, i.e., 
${\ovl{g}=e^{2u}g}$ with $\mathrm{grad}^N u=0$, on $M$. 

\begin{thm}\label{thm3}
Let $M^m \subset \wit{M}^{m+n}$ be a compact spin submanifold of a Riemannian spin manifold 
$(\wit{M},g)$.  Consider a non-trivial spinor field $\psi \in \gG(\bS)$ such that $D_H\psi=\gl\psi$. For any regular conformal change of the metric
$\ovl{g}=e^{2u}g$ on $\wit{M}$, assume that $$\ovl{R}e^{2u}+\kappa_1+4|Q^{\,\psi}|^2> ||H||^2 > 0$$ on $M_\psi$.
Then one has 
\begin{equation}\label{ethm3}
\gl^2 \geq \frac{1}{4}\; \inf _{M_{\psi}} \left( \sqrt {\ovl{R}e^{2u}+\kappa_1+4|Q^{\,\psi}|^2} -||H|| \right) ^2. 
\end{equation} 
\end{thm}

\begin{proof}
For $\psi \in \gG(\bS)$ an eigenspinor of $D_H$ with eigenvalue $\gl$, let 
$\ovl{\gp}:= e ^{-\frac{n-1}{2} u}\, \overline \psi$. Then (\ref{DHbar}) gives $\overline{D}_{H} \;\overline \gp = \gl\; e ^{-u}\; \overline  \gp\,.$
Recall that
\begin{eqnarray*}
\overline \gn _i \overline \psi = \overline {\gn _i \psi} -\frac{1}{2} \;\overline 
{e_i \cdot \mathrm{d} u \cdot \psi} -\frac{1}{2} e_i(u) \; \overline \psi, 
\end{eqnarray*}
and $\ovl{e_i}=e^{-u}e_i.$ Now, it is straightforward to get  
$
\overline Q^{\,\overline \gp}_{\bar i\, \bar j}= e ^{-u}Q ^{\,\psi}_{ ij}$, hence,
\begin{equation}\label{QQ}
|\overline Q^{\,\overline \gp}|^2=e^{-2u}|Q ^{\,\psi}|^2\;.
\end{equation}
Equation (\ref{Qpres}), which is also true on $(\wit{M},\ovl{g})$, applied to $\ovl{\gp}$ 
yields
\begin{align*}
\int_M& |\ovl{\gn}^{\ovl{Q}} \ovl{\gp} |^2v_{\ovl{g}}=\\&\int_M
(1+mq^2)\Big[(\gl e^{-u})^2-\frac{1}{4}\Big(\frac{\ovl{R}+\ovl{R}^N_{\ovl{\gp}}+4|\ovl{Q}^{\,\ovl{\gp}}| 
^2}{(1+mq^2)}-\frac{||\wit{H}||_{\ovl{g}}^2}{mq^2}\Big )\Big] |\ovl{\gp} |_{\ovl{g}}^2 \;v_{\ovl{g}}
\\-&\frac{1}{4}\int_M(1+mq^2)\Big[ \frac{2}{mq(1+mq^2)}(
||\wit{H}||_{\ovl{g}}^2-\frac{<\wit{H}\;\ovl{\cdot}\;\ovl{\gp},\ovl{\go_\perp}\;\ovl{\cdot}\;\ovl{\gp}>_{\ovl{g}}^2}{|\gp|_{\ovl{g}}^4})\Big]
|\ovl{\gp} |_{\ovl{g}}^2\; v_{\ovl{g}}\;.\end{align*}Since
$\wit{H}=e^{-u}\ovl{H}$, and
$\ovl{R}^N_{\ovl{\gp}}=e^{-2u}R^N_{\psi}$, we have
\begin{align*}\int_M& |\ovl{\gn}^{\ovl{Q}} \ovl{\gp} |^2v_{\ovl{g}}=\nonumber\\&\int_M
(1+mq^2) e^{-2u}\Big[\gl^2-\frac{1}{4}\Big(\frac{\ovl{R}e^{2u}+R^N_\psi+4|Q^{\,\psi}| 
^2}{(1+mq^2)}-\frac{||H||^2}{mq^2}\Big )\Big] |\ovl{\gp} |^2 v_{\ovl{g}}\nonumber\\-&\frac{1}{4}\int_M(1+mq^2)e^{-2u}\Big[ \frac{2}{mq(1+mq^2)}( ||H||^2-\frac{<H\cdot\psi,\go_\perp\cdot\psi>^2}{|\psi|^4})\Big] |\ovl{\gp} |^2 v_{\ovl{g}}\;. 
\end{align*}
As in the proof of Theorem \ref{thm2}, we finally take
$$q=\sqrt{\frac{||H||}{m(\sqrt{\ovl{R}e^{2u}+\kappa_1+4|Q^{\,\psi}|^2}-||H||)}}$$
and use the Cauchy-Schwarz inequality (\ref{CS}).
\end{proof}

If the hypothesis in Theorem \ref{thm3} is satisfied by an eigenfunction
$u_1$ associated with the first eigenvalue $\mu_1$ of the Yamabe
operator, then one has:

\begin{cor}
Under the same conditions as in Theorem \ref{thm3}, assume that $m\geq 3$ and
$\mu_1+\kappa_1+4|Q^\psi|^2>||H||^2>0$ on $M_\psi$, then
\begin{equation*}
\gl^2 \geq \frac{1}{4}\; \inf _{M_{\psi}} \left( \sqrt {\mu_1+\kappa_1+4|Q^{\,\psi}|^2} -||H|| \right) ^2. 
\end{equation*}
\end{cor} 
\begin{cor}
Under the same conditions as in Theorem \ref{thm3}, if $M$ is a
compact surface of genus zero and
$\frac{8\pi}{\mathrm{Area}(M)}+\kappa_1+4|Q^\psi|^2>||H||^2>0$ on $M_\psi$, then
\begin{equation*}
\gl^2 \geq \frac{1}{4}\; \inf _{M_{\psi}} \left( \sqrt {\frac{8\pi}{\mathrm{Area}(M)}+\kappa_1+4|Q^{\,\psi}|^2} -||H|| \right) ^2. 
\end{equation*}
\end{cor}

Now suppose that equality holds in (\ref{ethm3}). Then 
\begin{align*}
&\ovl{\gn}^{\ovl{Q}}\ovl{\gp}=0\;,\qquad
\go_\perp\cdot\psi=\varepsilon \frac{H}{||H||}\cdot\psi\quad \textrm{ where }\varepsilon\in\{\pm 1\}\;,\\
&|\gl|= \frac{1}{2}\left( \sqrt {\ovl{R}e^{2u}+\kappa_1+4|Q^{\,\psi}|^2} -||H|| \right)\quad\mathrm{and}\quad
\mathcal{R}^N\psi=\kappa_1\psi\;.
\end{align*}
Using (\ref{chang}) and (\ref{QQ}), it follows $\varepsilon=\mathrm{sgn}(\gl)$ and
\begin{equation}
\gn_i\psi^\star=\frac{1}{2}e_{i}\underset{M}\cdot \mathrm{d}u\underset{M}\cdot\psi^\star+\frac{m}{2}\mathrm{d}u(e_{i})\psi^\star-\sum_j 
Q^\psi_{ij}\,e_{j}\underset{M}\cdot\psi^\star\label{wem}
\end{equation} with $\mathrm{d}u=\frac{2\mathrm{d}(\mathrm{ln} (|\psi|))}{m-1}$. Non-trivial spinor fields satisfying equation (\ref{wem}) will be
naturally {\it called} twisted WEM-spinors (compare with \cite{Mor1}). 

\section{Final Remark}   

In this section, we show that the normal bundle of the submanifold can be replaced by an 
auxiliary arbitrary vector bundle on the submanifold. Thus, all the preceding computations could be done in an intrinsic way to obtain
results for a twisted Dirac operator on the manifold.

Let $(M^m,g)$ be a compact Riemannian spin manifold. Let $N\rightarrow M$ be a Riemannian vector bundle of rank $n$ over $M$. Suppose that $N$ is endowed with a 
metric connection $\gn^N$ and a spin structure. Let $\gS M$ (resp. $\gS N$) be the spinor 
bundle of $M$ (resp. $N$). Set $$\gS:=\gS M\otimes\gS N.$$ Recall that Clifford 
multiplication on $\gG(\gS)$ by a tangent vector field $X$ is given by: $$\forall 
\psi\in\gG(\gS),\quad X\cdot\psi=(\rho_M(X)\otimes\Id_{\gS N})(\psi).$$
Define the tensor-product connection $\gn$ on $\gG(\gS)$ by $$\gn=\gn^{\gS M}\otimes\Id_{\gS 
N}+\Id_{\gS M}\otimes\gn^{\gS N},$$ where $\gn^{\gS M}$ and $\gn^{\gS N}$ are the induced 
connections on $\gG(\gS M)$ and $\gG(\gS N)$ respectively. Let $D_M^{\gS N}$ be the twisted 
Dirac operator given by $$D_M^{\gS N}=\sum_{i=1}^{m} e_i\cdot\gn_i.$$ For any smooth real function 
$f$ on $M$, define the modified twisted Dirac operator by
$$D_f=D_M^{\gS N}-\frac{f}{2}.$$ For $\gl \in \rl$, consider the following modified connections 
\begin{eqnarray*}
\wih{\gn} _i ^\gl&=&\gn _i +\frac{(1-q)f}{2(1-mq)}\, e_i \cdot + q\gl e_i \cdot \\
\wih{\gn} _i ^Q&=&\gn _i - \frac{f}{2mq}\, e_i\cdot +(-1)^{n+1}q\gl e_i \cdot +\,\sum_j 
Q^{\,\psi}_{ij}\;e_j\cdot
\end{eqnarray*} 
where $Q^{\,\psi}$ is now the intrinsic Energy-Momentum tensor associated with $\psi$.

Note that these connections can be obtained from those defined in section~3, assuming 
that $$H\cdot\psi=f\go_\perp\cdot\psi.$$ In fact, this is the only way to give an intrinsic 
meaning to the modified connection used before. Then the same computations as in the proofs 
of Theorem \ref{thm1}, \ref{thm2} and \ref{thm3}, lead to the following assertions:

Let $(M^m,g)$ be a compact Riemannian spin manifold with $N\rightarrow M$ an auxiliary oriented Riemannian spin vector bundle of rank $n$.  
Let $\psi \in \gG(\gS)$ be an eigenspinor for the modified twisted Dirac 
operator $D_f$, associated with the eigenvalue $\gl$. Then,

\begin{prop}\label{gthm1}
Assume that $m\geq 2$ and $m 
(R+\kappa_1)>(m-1) f^2 > 0$ on $M_\psi$.
Then one has 
\begin{equation*}
\gl^2 \geq \frac{1}{4}\; \inf _{M_{\psi}} \left( \sqrt {\frac{m}{m-1}(R+\kappa_1)} -|f| \right) ^2. 
\end{equation*} 
If equality holds, $(M^m,g)$ admits a twisted Killing spinor. 
\end{prop}

Following the proof of Theorem \ref{thm3}, we can extend the previous theorem by performing a 
conformal change of the metric on $M$. For the limiting case, just note that 
${Q^\psi=\frac{1}{m}\mathrm{tr}(Q^\psi)g}$.

\begin{prop}\label{gthm2}
Assume that $m\geq 2$ and $m(\ovl{R}e^{2u}+\kappa_1)>(m-1) f^2 > 0$ on $M_\psi$ for any conformal change of the metric $\ovl{g}=e^{2u}g$ on $M$.
Then one has 
\begin{equation*}
\gl^2 \geq \frac{1}{4}\; \inf _{M_{\psi}} \left( \sqrt {\frac{m}{m-1}(\ovl{R}e^{2u}+\kappa_1)} -|f| \right) ^2. 
\end{equation*} 
If equality holds, $(M^m,\ovl{g})$ admits a twisted WEM-spinor, with $Q^{\,\psi}=\frac{\mu}{m}g$, where $$\mu^2=\quart\frac{m}{m-1}(\ovl{R}e^{2u}+\kappa_1)\;.$$ 
\end{prop}

\begin{prop}\label{gthm3} Assume that $R+\kappa_1+4|Q^{\,\psi}|^2> f^2 > 0$ on $M_\psi$.
Then one has 
\begin{equation*}
\gl^2 \geq \frac{1}{4}\; \inf _{M_{\psi}} \left( \sqrt {R+\kappa_1+4|Q^{\,\psi}|^2} -|f| \right) ^2. 
\end{equation*} 
If equality holds, $(M^m,g)$ admits a twisted EM-spinor. \end{prop} 

\begin{prop}\label{gthm4} Assume that $\ovl{R}e^{2u}+\kappa_1+4|Q^{\,\psi}|^2> f^2 > 0$ on $M_\psi$ for any conformal change of the metric $\ovl{g}=e^{2u}g$ on $M$.
Then one has 
\begin{equation*}
\gl^2 \geq \frac{1}{4}\; \inf _{M_{\psi}} \left( \sqrt {\ovl{R}e^{2u}+\kappa_1+4|Q^{\,\psi}|^2} -|f| \right) ^2. 
\end{equation*} 
If equality holds, $(M^m,g)$ admits a twisted WEM-spinor. 
\end{prop}

\begin{rem}
Assuming the normal curvature tensor is zero and $f$ is constant, then the necessary conditions for the equality cases in Theorem \ref{gthm1}, \ref{gthm2}, \ref{gthm3} and \ref{gthm4} become sufficient conditions. Moreover, when $m$ is odd, the considered Dirac operator may have to be defined with the opposite Clifford multiplication according to the sign of $f$.
\end{rem}  

\providecommand{\bysame}{\leavevmode\hbox to3em{\hrulefill}\thinspace}

\end{document}